\documentclass[10pt]{amsart}
\usepackage{enumerate,amsmath,amssymb,latexsym,
amsfonts, amsthm, amscd}

\theoremstyle{plain}

\nonumber

\begin{document}

\title {Cayley-Sylvester invariants and the Hamilton equations}

\date{}

\author[P.L. Robinson]{P.L. Robinson}

\address{Department of Mathematics \\ University of Florida \\ Gainesville FL 32611  USA }

\email[]{paulr@ufl.edu}

 \keywords{}

\begin{abstract}

We note implications of the Cayley-Sylvester theory of invariants and covariants for the Hamilton equations generated by cubic and quartic Hamiltonian functions. 

\end{abstract}

\maketitle

\bigbreak

\section*{} 

Let $u, v$ be a symplectic basis for the symplectic plane, in which an arbitrary point thereby decomposes as $z = p u + q v$ with $p$ and $q$ as linear symplectic coordinates. The Hamilton equations of motion generated by a Hamiltonian $\psi$ then read 
$$- \dot{p} = \psi_q, \; \; \dot{q} = \psi_p$$
and lead to the second-order equations 
$$\ddot{p} = \psi_{p q} \psi_q - \psi_p \psi_{qq}, \; \; \ddot{q} = \psi_{q p} \psi_p - \psi_q \psi_{pp}.$$
In what follows, we propose to consider an arbitrary solution curve to these Hamilton equations, the vector $z$ and its symplectic coordinates $p, q$ being functions of time. 

\medbreak 

When $\psi$ is the homogeneous cubic polynomial given by 
$$ 3 \psi = a p^3 + 3 b p^2 q + 3 c p q^2 + d q^4$$
the foregoing second-order equations in vector form read 
$$\ddot{z} = 2 F z$$
where $F$ is the quadratic in $p$ and $q$ given by 
$$F = (b^2 - ac) p^2 + (b c - a d) p q + (c^2 - b d) q^2.$$
Further differentiation along the solution curve reveals that this scalar function $F$ satisfies the second-order equation 
$$\ddot{F} = 6 F^2$$
and the first-order equation 
$$\dot{F}^2 = 4 F^3 - g_3$$
where $g_3$ is a constant of the motion. This situation is discussed quite thoroughly in [R]. 

\medbreak 

The classical theory of invariants and covariants developed largely by Cayley and Sylvester has a direct bearing on these matters; for this classical theory, we take [S] as our reference and import its notation. Associated to the homogeneous cubic 
$$U = a p^3 + 3 b p^2 q + 3 c p q^2 + d q^4$$
are the following: an invariant, its discriminant 
$$D = a^2 d^2 - 3 b^2 c^2 + 4 a c^3 + 4 b^3 d - 6 a b c d;$$
a quadratic covariant, its Hessian 
$$H = (ac - b^2) p^2 + (ad - bc) pq + (bd - c^2) q^2;$$
and a cubic (Jacobian) covariant 
$$J = (2 b^3 + a^2 d - 3 abc) p^3 + 3 (abd  +  b^2 c - 2 a c^2) p^2 q + 3 (2 b^2 d - b c^2 - acd) p q^2 + (3 bcd - a d^2 - 2 c^3) q^3.$$
See [S] page 183. These are related by an identity due to Cayley: 
$$J^2 = - 4 H^3 + D U^2$$
for which see [S] pages 186-187. 

\medbreak 

Now, if $U = 3 \psi$ then it is evident that $H = - F$ and a direct computation shows that $J = - \dot{F}$. Consequently, the Cayley identity $J^2 = - 4 H^3 + D U^2$ becomes $\dot{F}^2 = 4 F^3 + D (3 \psi)^2$ or 
$$\dot{F}^2 = 4 F^3 - g_3$$
where 
$$g_3 = - D (3 \psi)^2.$$

\medbreak 

This observation has prompted us to look further; here, we shall present only a discussion of quartic Hamiltonians. Thus, let $\psi$ be the homogeneous quartic polynomial given by 
$$4 \psi = a p^4 + 4 b p^3 q + 6 c p^2 q^2 + 4 d p q^3 + e q^4.$$
The Hamilton equations for this Hamiltonian lead to the second-order equations in vector form 
$$\ddot{z} = 3 G z = \frac{3}{4} F z$$
where the scalar function $G$ is given by 
$$G = (b^2 - ac) p^4 + 2(bc - ad) p^3 q + (3 c^2 - 2 bd - ae) p^2 q^2 + 2 (cd - be) p q^3 + (d^2 - ce) q^4.$$
Further differentiation reveals that $F = 4 G$ satisfies the first-order equation 
$$\dot{F}^2 = 4 F^3 - g_2 F - g_3$$
where $g_2$ and $g_3$ are constants of the motion. We omit the somewhat lengthy justification of this, as the outcome reproduces the following items of classical invariant theory. 

\medbreak 

Again we switch to the notation of [S]. Associated to the homogeneous quartic 
$$U = a p^4 + 4 b p^3 q + 6 c p^2 q^2 + 4 d p q^3 + e q^4$$ 
are invariants due to Cayley
$$S = ae - 4 bd + 3 c^2$$ 
and to Boole 
$$T = ace + 2 bcd - a d^2 - b^2 e - c^3$$
along with two covariants: a quartic 
$$H = (ac - b^2) p^4 + 2 (ad - bc) p^3 q + (ae + 2bd - 3 c^2) p^2 q^2 + 2(be - cd) p q^3 + (ce - d^2) q^4$$
and a sextic $J$ which we refrain from displaying in full. For details, see [S] pages 189 and 192. These invariants and covariants are related by the identity
$$J^2 = - 4 H^3 + S U^2 H - T U^3$$
for which see [S] pages 195 and 200. 

\medbreak 

Now, if we take $U = 4 \psi$ then it turns out that $H$ becomes a multiple of $- F$ and $J$ a multiple of $ - \dot{F}$; the identity relating these covariants to $S$ and $T$ then assumes the form
$$\dot{F}^2 = 4 F^3 - g_2 F - g_3$$
where in fact 
$$g_2 = S (16 \psi)^2,$$
$$g_3 = T (16 \psi)^3.$$

\medbreak 

Before closing, we offer some additional comments specific to cubics and quartics. 

\medbreak 

{\it Cubics}: Barring obvious exceptions, the solution $F$ to $\dot{F}^2 = 4 F^3 - g_3$ is a (shifted) Weierstrass $\wp$ function; as $g_2$ is zero, this is the equianharmonic case corresponding to a triangular lattice. The second-order equation $\ddot{z} = 2 F z$ is then of Lam\'e type with $n = 1$ and may be solved as such; see [F] page 285. 

\medbreak 

{\it Quartics}: The typical solution to $\dot{F}^2 = 4 F^3 - g_2 F - g_3$ is also a (shifted) Weierstrass $\wp$ function; the nature of the corresponding lattice is determined by
$$g_2^3 - 27 g_3^2 = (S^3 - 27 T^2) (16 \psi)^6$$
where $S^3 - 27 T^2$ is the discriminant of the quartic; see [S] page 191. Here, the second--order equation $\ddot{z} = 2 F z$ is of Lam\'e type with $n = 1/2$; see [WW] page 578. 

\medbreak 

Finally, we remark that a symplectic connexion to the invariant theory of binary quantics in general is manifest starting from Article 140 of [S].

\bigbreak

\begin{center} 
{\small R}{\footnotesize EFERENCES}
\end{center} 
\medbreak 

[F] A.R. Forsyth, {\it Theory of Functions of a Complex Variable}, Cambridge, First Edition (1893). 

[R]  P.L. Robinson, {\it Cubic Hamiltonians}, arXiv 1505.00671 (2015).

\medbreak

[S] G. Salmon, {\it Lessons Introductory to the Modern Higher Algebra}, Dublin, First Edition (1885). 

\medbreak 

[WW] E.T. Whittaker and G.N. Watson, {\it A Course of Modern Analysis}, Cambridge, Fourth Edition (1950). 

\medbreak

\end{document}